\documentclass[a4paper, twoside, 10pt]{article}
\usepackage{amsmath,amsthm,amsfonts,latexsym,amscd,amssymb, enumerate}
\usepackage{hyperref}
\numberwithin{equation}{section}
\newtheorem{theorem}{Theorem}[section]
\newtheorem{lemma}{Lemma}[section]
\newtheorem{corollary}{Corollary}[section]
\newtheorem{proposition}{Proposition}[section]
\newtheorem{definition}{Definition}[section]
\newtheorem{example}{Example}[section]
\theoremstyle{remark}

\date{}
\title{\textbf{Conformal Warped Product Submersion}}
\author{Harmandeep Kaur, Abhishek Pandey and Gauree Shanker\thanks{corresponding author, Email: grshnkr2007@gmail.com}}

\begin{document}
	\maketitle
\begin{center}
	\textbf{Abstract}
\end{center}
In this paper, the concept of Riemannian warped product submersion is generalized to the conformal case. We introduce the notion of conformal warped product submersion. It is a submersion between warped product manifolds that preserves angles between the horizontal vectors. The fundamental tensors of submersion are derived for conformal warped product submersion.

\section{Introduction}
In order to compare the geometric structures between two manifolds, we need suitable types of maps between Riemannian manifolds. Given two manifolds, the maps are known as submersions if the rank of a differential map is equal to the dimension of the target manifold and immersions if the rank of a differential map is equal to the dimension of source manifold. Moreover, if these maps are isometry between manifolds, then the immersion is called isometric immersion  and the submersion
is called Riemannian submersion. Riemannian submersions were introduced in the sixties by Gray \cite{J13} and O’Neill \cite{J8}. Riemannian submersion is a tool to study the geometry of a Riemannian manifold with an additional structure in terms of certain components, that is, the fibers and the base space. Riemannian submersions are related to physics and have their applications in the Yang–Mills theory \cite{J6,j1}, Kaluza–Klein theory \cite{J7,j2}, supergravity and superstring theories \cite{j3,j4}. The projection of a Riemannian product manifold on one of its factors is a trivial example of Riemannian submersion.\\
The class of warped product manifolds has shown itself to be rich, both wide and diverse, playing important roles in differential geometry as well as in physics. To illustrate, Bishop and O’Neill introduced warped products in \cite{J1} as means to construct a large class of complete Riemannian manifolds with negative curvature. The notion of warped product manifolds is one of the most fruitful generalizations of Riemannian products. Such notion plays very important roles in differential geometry as well as in physics, especially in general relativity. Schwarzschild and Robertson Walker cosmological models are well-known examples of warped product manifolds \cite{J20}. Warped product manifolds have been studied for a long period of time. In contrast, the study of warped product submanifolds was only initiated around the beginning of this century in a series of articles \cite{j51, j36, j53, j54}.\\
The study of maps between Riemannian warped product manifolds is an active research field. The theory of warped product immersion has been studied extensively so far. On the other hand, the study of Riemannian warped product submersions is an emerging area of research. The notion of Riemannian warped product submersion was introduced by I. K. Erken and C. Murathan in $2021$ \cite{J11}. Recently, I. K. Erken et al. extended the study on Riemannian warped product submersion and gave the curvature properties of such submersions in \cite{j20}.\\
Conformal submersion is a generalization of Riemannian submersion. It is a submersion between Riemannian manifolds having conformal mapping between horizontal vectors. Although conformal maps do not preserve the distance between points contrary to isometries, they preserve angles between vector fields. This property enables one to transfer certain properties of a manifold
to another manifold by deforming such properties. Conformal submersions  were introduced independently by Fuglede \cite{J15} and Ishihara \cite{J16} which is useful for the characterization of harmonic morphisms \cite{j22} and have applications in medical imaging and computer graphics. In \cite{j28}, Ornea obtained the fundamental equations of such submersions. The curvature relations for conformal submersions were given in \cite{j22, j21}. Further, conformal submersion was studied by many authors \cite{j33, j35, j31, j32, j34}. R. Tojeiro studied conformal immersions of warped product manifolds in \cite{j18}. So, it is interesting to study conformal warped product submersion. In this paper we introduce the notion of conformal warped product submersion which is the generalization of Riemannian warped product submersion.

\section{Preliminaries}
In this section, we are going to recall the foundational concepts required to understand the notion of conformal warped product submersion.

\subsection{Warped Product Manifolds}
\begin{definition}Let $(M_1^{m_1},g _{M_1})$ and $(M_2^{m_2}, g_{M_2})$ be two Riemannian manifolds and let f be a positive smooth function on $M_1$. Then the warped product $M = M_1 \times _{f}  M_2$ of $M_1$ and $M_2$ is the product manifold $M_1 \times M_2$ endowed with the metric $g_M$ defined as  
	\begin{equation}
		g_{M}(X ,Y) = g_{M_1}(\pi _{1\ast }(X), \pi _{1 \ast }(Y) ) + f ^{2}(\pi _{1}(x, y) )g _{M _2}(\pi _{2\ast }(X), \pi _{2 \ast }(Y) ),
	\end{equation} 
	where X, Y are vector fields on $M_1 \times M_2$ and $\pi_1$, $\pi_2$ are projection mappings of $M$ onto $M_1$, $M_2$ respectively.
\end{definition}
The fibers $\{x\} \times  M_2 = \pi_{1} ^{-1}(x)$ and the leaves $ M_1\times\{y\}=\pi_{2}^{-1}(y)$ are Riemannian submanifolds of $M_1 \times_{f} M_2$.
Vectors tangent to leaves and those tangent to fibers are called horizontal and vertical respectively \cite{j19}. If v $\in \ T_{p}M_1$ and $q \in M_2$, then the lift $\overline{v}$ of $v$ to $(p, q)$ is the unique vector $T_{(p, q)}M$ such that $(\pi_{1})_ {\ast } (\overline{v}) = v$. For a vector field X $\in \ \Gamma (TM_1)$, the lift of X to M is the vector field $\overline{X}$ whose value at each (p, q) is the lift of $X_p$ to (p, q). The set of all horizontal lifts is denoted by $\mathcal{L}(M_1)$. Similarly, we denote the set of all vertical lifts by $\mathcal{L}(M_2)$. A vector field $ \overline{E}$ of $M_1 \times \ M_2$ can be written as $\overline{E}= \overline{X} + \overline{U}$ with $\overline{X} \in \mathcal{L}(M_1)$ and $\overline{U} \in \mathcal{L}(M_2)$. 
\begin{lemma} \cite{J20}
	Let $M=M_{1}\times_{f} M_{2}$ be a warped product manifold. For any $E_1,F_1\in \mathcal{L}(M_1)$ and $E_2,F_2\in \mathcal{L}(M_2)$,
	\begin{enumerate}
		\item $\nabla_{E_1}F_1$ is the lift of $\nabla_{E_1}^{1}F_1$,
		\item $\nabla_{E_1}E_2=\nabla_{E_2}E_1=(E_{1}(f)/f)E_2$,
		\item $nor(\nabla_{E_2}F_2)=-g_M(E_2,F_2)(D \ ln \ f)$,
		\item $tan(\nabla_{E_2}F_2)$ is the lift of $\nabla_{E_2}^{2}F_2$.
	\end{enumerate}
\end{lemma}
Here $\nabla,\nabla^1$ and $\nabla^2$ denote Riemannian connections on $M,M_1$ and $M_2$, respectively and $Df$ denotes the gradient of f.
\begin{corollary} \cite{J20}
	Let $M = M_1 \times_f M_2$ be a warped product manifold. Then
	the leaves $M_1 \times \{y\}$ and the fibers $\{x\} \times M_2$ are totally geodesic and totally umbilical, respectively.
\end{corollary}  

\subsection{Riemannian submersions}
\begin{definition}\cite{J8}
	Let $(M^{m},g _{M})$ and $(N^{n}, g_{N})$ be two Riemannian manifolds, , where $dim(M) = m,$ $dim(N) = n$ and $m > n$. A Riemannian submersion $F: M \to N$ is a surjective map
	of M onto N satisfying the following axioms:
	\begin{enumerate}
		\item [(i)] F has maximal rank.
		\item [(ii)] The differential $F_\ast$ preserves the lengths of horizontal vectors.
	\end{enumerate}
\end{definition}
For each $b\in N$, $F^{-1}(b)$  is a submanifold of $M$ of dimension $(m - n)$, called the fibers of $F$. A vector field on $M$ is called vertical if it is always tangent to the fiber and is called horizontal if it is always orthogonal to the fibers.
The integrable distribution of $F$ is defined by $\mathcal{V}_p=ker F_{*p}$ and  $\mathcal{V}_p$ is called the vertical distribution of submersion $F$, and the distribution $\mathcal{H}_p= (\mathcal{V}_p)^\perp$ which is a complementary and orthogonal distribution to $\mathcal{V}$ is called horizontal distribution. Thus, for every $p \in M$, $M$ has the following decomposition:
\begin{equation}
	T_pM=\mathcal{V}_p\oplus \mathcal{H}_p=\mathcal{V}_p\oplus\mathcal{V}_p^\perp.
\end{equation}
A vector field $X$ on $M$ is called basic if $X$ is horizontal and $F$ related to a vector field $X_*$ on $N$, i.e., $F_*X_p = X_{*F(p)}$ for all $p \in M$.

The geometry of Riemannian submersions is characterized by O’Neill’s tensors \cite{J8} $\mathcal{T}$ and $\mathcal{A}$ defined for vector fields $E$, $F$ on $M$ by
\begin{equation}
	\mathcal{A}_{E}F = \mathcal{H}\nabla_{\mathcal{H}E} \mathcal{V}F + \mathcal{V}\nabla_{\mathcal{H}E} \mathcal{H}F,
\end{equation}   
\begin{equation}
	\mathcal{T}_{E}F = \mathcal{H}\nabla_{\mathcal{V}E} \mathcal{V}F + \mathcal{V}\nabla_{\mathcal{V}E} \mathcal{H}F,
\end{equation}
where, $\nabla$ is the Levi-Civita connection of $g_M$ and $\mathcal{T}$ acts as the second fundamental form of all the fibers and the tensor $\mathcal{A}$ determines the integrability of horizontal distributions. A Riemannian submersion is called a Riemannian submersion with totally geodesic fiber if $\mathcal{T}$ vanishes identically.
A Riemannian submersion is called a Riemannian submersion with totally umbilical fibre if
\begin{equation}
	\mathcal{T}_UW=g_M(U,W)H
\end{equation}
$\forall$ U, W $\in \Gamma (\mathcal{V}).$

\subsection{Conformal Submersion}
Let $ (M,g_{M}) $ and $ (B, g_{B}) $  be Riemannian manifolds and $ F: M \rightarrow B $ be a smooth submersion, then F is called a conformal submersion, if there is a positive function $\lambda :M \rightarrow \mathbb{R}^{+} $ such that
\begin{equation}
	\lambda^{2} g_{M}(X, Y )=g_{B}(F_{\ast} X,F_{\ast}Y )
\end{equation} 
for  X, Y $\in \Gamma (( kerF_{\ast })^{\perp })$  and $\lambda$ is called dilation.
\begin{proposition}
	\cite{j21} Let $ \pi : (M^{m},g) \rightarrow (N^{n},h)$ be a conformal submersion with dilation $\lambda$ and X, Y be horizontal vectors, then
	\begin{equation}
		A_{X}Y = \frac{1}{2} \bigg\{  \mathcal{V} [X, Y] -\lambda^{2}g(X, Y) grad_{\mathcal{V}}\Big( \dfrac{1}{\lambda ^2}\Big) \bigg\}. 
	\end{equation}
\end{proposition}

\begin{example} Let F be a map defined by, $$F:R^4\rightarrow R^2; (x_1,x_2,x_3,x_4)\rightarrow (e^{x_3}sin\ x_4,e^{x_3}cos\ x_4),$$then it is a conformal submersion with $\lambda =e^{x_3}.$
\end{example}

\subsection{Riemannian warped product submersion}
\begin{definition} \cite{J11}
	Let $\phi _{i}$, i = 1, 2, be Riemannian submersions from $M_i$ to $N_{i}$. If $M = M_1 \times_{f} \ M_2$ and $N = N_1\times_{ \rho }\ N_2$ are Riemannian warped product manifolds, then the map 
	\begin{equation}
		\phi = \phi _{1} \times \phi _{2} : M = M_1 \times_{f} M_2 \rightarrow N = N_1 \times _{\rho }  N_2 
	\end{equation}
	given by $( \phi _{1} \times \phi _{2})(p_{1}, p_{2}) = (\phi _{1}(p_1), \phi _{2}(p_{2})) $ is a Riemannian submersion. This kind of Riemannian submersions are called Riemannian warped product submersions.
\end{definition}

A Riemannian warped product submersion  $ \phi =( \phi _{1}, \phi _{2}) : M = M_1 \times_{f} M_2 \rightarrow N = N_1 \times _{\rho }  N_2$ has $M_i$-minimal fibers if $H_i$ vanishes identically for $i=1,2$ and mixed totally geodesics fibers if its second fundamental form T satisfies $T(E, F)=0$ for any $E\in \Gamma (\mathcal{V}_{1})$ and $F\in\Gamma(\mathcal{V}_{2}).$

\section{Conformal warped product submersion}
In this section, we introduce the notion of conformal warped product submersion which is a generalization of Riemannian warped product submersion. It is a submersion between warped product manifolds that preserve the angles between the horizontal vectors. Further, the expressions for fundamental tensor of submersion are derived for conformal warped product submersion.

\begin{definition}
	Let $M_1, M_2$ be Riemannian manifolds and $\lambda_1, \lambda_2$ be smooth positive functions on $M_1$ and $M_2,$ repectively. If $\lambda$ is a smooth positive function on $M_1 \times \ M_2$ i.e.,
	\begin{equation}
		\lambda : M_1 \times \ M_2 \rightarrow R^+ 
	\end{equation}
	such that , $\lambda \lvert _{M_1} = \lambda _1$ and $ \lambda \lvert_{M_2} = \lambda _2,$ then  $\lambda$ is called the lift function of $\lambda_1$ and $\lambda_2.$
\end{definition}
\begin{proposition}
	Let $\phi _{i}: M_i \rightarrow N_i$, i = 1, 2, be conformal submersions with dilation $\lambda _i$,  from Riemannian manifolds $M_i$ to $N_{i}$. If $M = M_1 \times _{f} \ M_2$ and $N = N_1\times_{ \rho }\ N_2$ are Riemannian warped product manifolds, then the map 
	\begin{equation}
		\phi = \phi _{1} \times \phi _{2} : M = M_1 \times_{f} \ M_2 \rightarrow N = N_1 \times_{\rho } \ N_2 
	\end{equation}
	given by $( \phi _{1} \times \ \phi _{2})(p_{1}, p_{2}) = (\phi _{1}(p_1), \phi _{2}(p_{2})) $ is a conformal submersion  with dilation $\lambda$, where $\lambda$ is the lift of $\lambda_1$ and $\lambda_2$ to $M_1 \times \ M_2.$ We call this  submersion as conformal warped product submersion. 
	\begin{proof}
		Let $X_1, Y_1$ and $X_2, Y_2$ be horizontal vector fields on $M_1$ and $M_2$ repectively. As $\phi_i :M_i \to N_i$ is conformal submersion for $i=1, 2$. So, from eq.$(2.6)$ we have
		\begin{equation}
			g_{N_i}(\phi_{i\ast}X_i, \phi_{i\ast}Y_i)= \lambda_{i}^2g_{M_i}(X_i, Y_i).  
		\end{equation}
		Since, $\phi_i$ is a submersion from $M_i$ to $N_i$ for $i=1, 2,$
		the map $\phi : M = M_1 \times_{f} \ M_2 \rightarrow N = N_1 \times_{\rho} \ N_2 $ is a submersion.
		\\
		The tangent space of M has the following decomposition for $p = (p_1 , p_2) \in M,$ where $p_1 \in M_1$ and $p_2 \in M_2,$
		\begin{equation}
			T _{(p_1 , p_2)}\big(M_1 \times M_2\big) = T _{(p_1 , p_2)}\big(M_1 \times \{ p_2 \} \big) \oplus T_{(p_1,p_2)}\big( \{p_1 \} \times  M_2\big),
		\end{equation}
		\begin{equation}
			T _{(p_1 , p_2)}\big(M_1 \times  M_2\big) = \mathcal{H} _{(p_1 , p_2)}  \oplus  \mathcal{V} _{(p_1 , p_2)},
		\end{equation}
		where, $\mathcal{H}, \mathcal{V}$ denotes the horizontal and vertical subspace of M respectively.\\
		Also, $ker( \phi _{1} \times \phi _{2}) _{\ast }  = ker( \phi _{1 \ast}) \times ker( \phi _{2 \ast })$\\
		Using $(3.4)$ and $(3.5),$ we have \\
		$T _{(p_1 , p_2)}\big(M_1 \times \{ p_2 \} \big) = \big((\mathcal{H}_1)_{p_1} \times \{ p_2 \} \big) \oplus \big((\mathcal{V}_1)_{p_1} \times \{ p_2 \} \big), $\\
		$T _{(p_1 , p_2)}\big( \{p_1 \} \times  M_2\big) = \big( \{p_1 \} \times (\mathcal{H}_2)_{p_2} \big) \oplus \big( \{p_1 \} \times (\mathcal{V}_2)_{p_2} \big)$.\\
		Hence, we get the decomposition of vertical and horizontal subspace of M as,
		\begin{equation*}
			\begin{split}
				\mathcal{V} _{(p_1 , p_2)} &= \big((\mathcal{V}_1)_{p_1} \times \{ p_2 \} \big) \oplus \big( \{p_1 \} \times (\mathcal{V}_2)_{p_2} \big),\\
				\mathcal{H} _{(p_1 , p_2)} &= \big((\mathcal{H}_1)_{p_1} \times \{ p_2 \} \big) \oplus \big( \{p_1 \} \times (\mathcal{H}_2)_{p_2} \big).
			\end{split}
		\end{equation*}
		For a horizontal vector field $X_{i}^{\mathcal{H}} \in \Gamma(\mathcal{H}_i),$ the lift of $X_{i}^{\mathcal{H}}$ to M is the vector field $(\overline{X_{i}^{\mathcal{H}}}) =(\overline{X_{i}})^\mathcal{H}$. Similarly, for a vertical vector field  $X_{i}^{\mathcal{V}} \in \Gamma(\mathcal{V}_i),$ the lift of $X_{i}^{\mathcal{V}}$ to M is the vector field $( \overline{X_{i}^{\mathcal{V}}} ) =(\overline{X_{i}})^\mathcal{V}.$\\
		For instance, both the vector field on $M_i$ and its lift to M will be denoted by the same notation, the meaning will be clear from the context.\\
		Now, in order to show that the submersion, $\phi: M \to N$ is a conformal submersion, we proceed as follows: 
		\begin{equation*}
			\begin{split}
				g_{N}(\phi_{\ast}(X_1,X_2), \phi_{\ast}(Y_1,Y_2))&=g_{N_1}(\phi_{1\ast}(X_1), \phi_{1\ast}(Y_1))+ \rho^2(\phi_1(p_1))g_{N_2}(\phi_{2\ast}(X_2), \phi_{2\ast}(Y_2))\\
				& =\lambda_{1}^2(p_1)g_{M_1}(X_1, Y_1)+ f^2(p_1).\lambda_{2}^2(p_2)g_{M_2}(X_2, Y_2)\\
				&=\lambda_{1}^2(\pi_1(p))g_{M_1}(X_1, Y_1)+ f^2(p_1).\lambda_{2}^2(\pi_2 (p))g_{M_2}(X_2, Y_2)\\
				&=(\lambda_{1}\circ \pi_1)^2g_{M_1}(X_1, Y_1)+ f^2(p_1).(\lambda_{2}\circ \pi_2)^2g_{M_2}(X_2, Y_2)\\
				& =\lambda^2(g_{M_1}(X_1, Y_1)+ f^2g_{M_2}(X_2, Y_2))\\
				&=\lambda^2g_M((X_1, X_2), (Y_1, Y_2)),
			\end{split}
		\end{equation*}
		where $\lambda$ is the lift function of $\lambda_1$ and $\lambda_2.$\\
		Therefore, $\phi : M = M_1 \times_{f} \ M_2 \rightarrow N = N_1 \times_{\rho } \ N_2 $ is a conformal warped product submersion with dilation $\lambda$.
	\end{proof}
\end{proposition}
\begin{corollary}
	If the dilation $\lambda \equiv 1$, then the conformal warped product submersion is a Riemannian warped product submersion. Thus, Riemannian warped product submersion is a particular case of conformal warped product submersion.
\end{corollary}

\begin{corollary}
	If $\phi$ is the conformal warped product submersion with dilation $\lambda = e^{- \sigma} ; \sigma \in C^{\infty}(M),$ then the metric $G_M$ given by $G_M = e^{2\sigma} g_M$ is the unique metric on $M = M_1 \times_{f} M_2$ conformal with $g_M$ with the property that
	\begin{equation}
		\varphi : (M, G_M) \to (N, g_N),
	\end{equation}
	where $\varphi$ is a Riemannian submersion, defined by $\varphi(p) = \phi(p)$, for $p\in M.$ 
	\begin{proof}
		Let $X_1, Y_1$ and $X_2, Y_2$ be horizontal vector fields on $M_1$ and $M_2$ repectively. Since, $\phi$ is the conformal warped submersion with dilation $\lambda.$ By Prop. $3.1,$ we have 
		\begin{equation*}
			\begin{split}
				g_{N}(\phi_{\ast}(X_1,X_2), \phi_{\ast}(Y_1,Y_2)) &= \lambda^2g_M\big((X_1, X_2), (Y_1, Y_2)\big)\\
				&= e^{-2\sigma}g_M\big((X_1, X_2), (Y_1, Y_2)\big)\\
				&= e^{-2\sigma}e^{2\sigma}G_M\big((X_1, X_2), (Y_1, Y_2)\big)\\
				&= G_M\big((X_1, X_2), (Y_1, Y_2)\big).
			\end{split}
		\end{equation*}
		Hence, we get the assertion.   
	\end{proof}
\end{corollary}

\begin{theorem}
	Let $ \phi =( \phi _{1}, \phi _{2}) : M = M_1 \times_{f} M_2 \rightarrow N = N_1 \times _{\rho }  N_2$ be a  conformal  warped product submersion with dialtion $\lambda$ between two Riemannian warped product manifolds. If  $X_i,Y_i \in \ \Gamma (\mathcal{H}_{i}), i=1,2$ on $M_i,$ then  we have
	\begin{enumerate}
		\item  $A(X_1, Y_1) = A_1(X_1, Y_1)$ = $\dfrac{1}{2} \Big\{  \mathcal{V} [X_1, Y_1] -\lambda_{1}^{2}g_{M_1}(X_1, Y_1 ) grad_{\mathcal{V}}\Big( \dfrac{1}{\lambda_{1} ^2}\Big) \Big\}, $
		\item $A(X_2, Y_2) = \dfrac{1}{2} \Big\{  A_2(X_2, Y_2) - A_2(Y_2, X_2) -  \lambda_{2}^{2}g_{M_2}(X_2, Y_2 ) grad_{\mathcal{V}}\Big( \dfrac{f^2}{\lambda_{1} ^2}\Big) \Big\}. $
	\end{enumerate}
	\begin{proof}
		Extend $X_i, Y_i$ to basic vector fields.\\
		\hspace{0.3cm} In view of $(2.7)$, we have
		\begin{equation}
			\begin{split}
				A(X_1, Y_1)&= \frac{1}{2} \{  \mathcal{V} [X_1, Y_1] - grad_{\mathcal{V}}( g_M(X_1, Y_1) ) \} \\
				&= \frac{1}{2} \{  \mathcal{V} [X_1, Y_1] - grad_{\mathcal{V}}( g_{M_1}(X_1, Y_1) ) \}.
			\end{split}
		\end{equation}
		Combining $(2.6)$ and $(3.7),$ we get,
		\begin{equation*}
			\begin{split}
				A(X_1, Y_1) &= \dfrac{1}{2} \Big\{  \mathcal{V} [X_1, Y_1] - grad_{\mathcal{V}}\Big(\dfrac{1}{\lambda _{1}^{2}}g_{N_1}(\phi_{1\ast}(X_1), \phi_{1\ast}(Y_1))\Big)\Big\}\\
				&= \dfrac{1}{2} \Big\{ \mathcal{V} [X_1, Y_1] - g_{N_1}\big(\phi_{1\ast}(X_1), \phi_{1\ast}(Y_1)\big)grad_{\mathcal{V}}\Big(\frac{1}{\lambda _{1}^{2}}\Big) \Big\}\\
				&= \dfrac{1}{2} \Big\{  \mathcal{V} [X_1, Y_1] - \lambda _{1}^{2}g_{M_1}\big(\phi_{1\ast}(X_1), \phi_{1\ast}(Y_1)\big)grad_{\mathcal{V}}\Big(\dfrac{1}{\lambda _{1}^{2}}\Big) \Big\}\\
				&= A_1(X_1, Y_1),
			\end{split}
		\end{equation*}
		which proves (1).\\
		For (2), we proceed as follows
		\begin{equation}
			\begin{split}
				A(X_2, Y_2) &= \dfrac{1}{2}\{ \mathcal{V} [X_2, Y_2] - grad_{\mathcal{V}}\big( g_{M}(X_2, Y_2) \big)\} \\
				&= \frac{1}{2} \{  \mathcal{V} [X_2, Y_2] - grad_{\mathcal{V}}\big( f^2g_{M_2}(X_2, Y_2) \big) \}.
			\end{split}
		\end{equation}
		Combining $(2.6)$ and $(3.8),$  we get
		\begin{equation*}
			A(X_2, Y_2) = \frac{1}{2} \Big\{  \mathcal{V} (\nabla_{X_2}Y_2 - \nabla_{Y_2}X_2) - grad_{\mathcal{V}}\Big(\dfrac{f^2}{\lambda _{2}^{2}}g_{N_2}(\phi_{2\ast}(X_2), \phi_{2\ast}(Y_2)) \Big)\Big\}.
		\end{equation*}
		Using lemma $(2.1)$ and above equation we obtain
		\begin{equation*}
			\begin{split}
				A(X_2, Y_2) &= \dfrac{1}{2} \Big\{  \mathcal{V} ( (\nabla_{X_2}^2Y_2 - g_M(X_2, Y_2)(grad \ ln \ f))\\
				&- (\nabla_{Y_2}^2X_2
				- g_M(Y_2, X_2)(grad \ ln \ f)))\\
				& -g_{N_2}(\phi_{2\ast}(X_2), \phi_{2\ast}(Y_2)) grad_{\mathcal{V}}\Big(\ \dfrac{f^2}{\lambda _{2}^{2}}\Big) \Big\}\\
				&= \dfrac{1}{2} \Big\{  \mathcal{V} (\nabla_{X_2}^2Y_2 - \nabla_{Y_2}^2X_2 ) - \lambda_{2}^2g_{M_2}(X_2, Y_2) grad_{\mathcal{V}}\Big(\ \dfrac{f^2}{\lambda _{2}^{2}}\Big) \Big\}\\
				&= \dfrac{1}{2} \Big\{ ( A_{2}(X_2, Y_2 ) - A_{2}(Y_2, X_2 ) - \lambda_{2}^2g_{M_2}(X_2, Y_2) grad_{\mathcal{V}}\Big(\ \dfrac{f^2}{\lambda _{2}^{2}}\Big) \Big\}, 
			\end{split}
		\end{equation*}
		which proves (2).
	\end{proof}
\end{theorem}
In our study, the concept of Riemannian warped product submersion is generalized to the conformal case which leads number of interesting and significant issues to investigate. 

\section*{Acknowledgments}
The first author is thankful to UGC for providing financial assistance in terms of JRF scholarship vide NTA Ref. No.: $201610070797$(CSIR-UGC NET June 2020). The third author is thankful to the Department of Science and Technology (DST) Government of India for providing financial assistance in terms of FIST project (TPN-69301) vide the letter with Ref. No.: (SR/FST/MS-1/2021/104).

\noindent 	H. Kaur, A. Pandey and G. Shanker\newline
Department of Mathematics and Statistics\newline
Central University of Punjab\newline
Bathinda, Punjab-151401, India.\newline
Email: harmandeepkaur1559@gmail.com, abhishekpandey950498@gmail.com, grshnkr2007@gmail.com
\end{document}